\documentclass[aps,pra,amssymb,amsmath,eqsecnum,nofootinbib]{revtex4}
\usepackage{graphicx}

\newtheorem{definition}{Definition}[section]
\newtheorem{theorem}{Theorem}[section]

\newtheorem{lemma}{Lemma}[section]

\newtheorem{axiom}{AXIOM}[section]

\newenvironment{hypothesis}{HP: \begin{center}} {\end{center}}
\newenvironment{thesis}{TH: \begin{center}} {\end{center}}
\newenvironment{proof}{\begin{center}PROOF: \end{center}} {$ \blacksquare $}

\begin{document}
\title{An interesting temporalization of G\"{o}del's ontological proof}
\author{Gavriel Segre}
\homepage{http://www.gavrielsegre.com}
\begin{abstract}
Recent theologies concerning God's death after Auschwitz are mathematically formalized through a suitable temporalization
of G\"{o}del's Ontological Proof.
\end{abstract}
\maketitle
\newpage
\tableofcontents
\newpage
\section{Acknowledgements}

First of all I would like to thank Piergiorgio Odifreddi for many stimulating discussions
concerning the philosophical and mathematical meaning of G\"{o}del's ontological proof.

 I would then like to thank Vittorio de Alfaro for his friendship and
his moral support, without which I would have already given up.

Then  I would like to thank Andrei Khrennikov and the whole
team at the International Center of Mathematical Modelling in
Physics and Cognitive Sciences of V\"{a}xj\"{o} for their very
generous informatics' support.

Of course nobody among the mentioned people has responsibilities
as to any (eventual) error contained in these pages.

Last but not least I dedicate this 25-April-paper to my first Teacher (in the meaning of Pirk\'{e} Avot) Gianni Jona-Lasinio.

\newpage
\section{Modal logic}
Let us introduce briefly Charles Lewis' formal systems of Modal
Logic \cite{Konyndyk-86}.

Introduced the \emph{necessity operator} $ \Box $ and the
\emph{possibility operator} $ \lozenge $ let us introduce the
following:

\begin{definition}
\end{definition}
\emph{T formal system:}

the formal system obtained adding to the propositional logic the
modal operators $ \Box $ and  $ \lozenge $ and the following
axioms:
\begin{equation}
  \Box \phi \; \rightarrow \; \phi
\end{equation}
\begin{equation}
    \phi \; \rightarrow \; \lozenge \phi
\end{equation}
\begin{equation}
    \Box \phi \; \leftrightarrow \; \neg \lozenge \neg \phi
\end{equation}
\begin{equation}
    \lozenge \phi \; \leftrightarrow \; \neg \Box \neg \phi
\end{equation}
\begin{equation}
    \Box ( \phi \wedge \psi) \; \leftrightarrow \; \Box \phi \wedge \Box \psi
\end{equation}
\begin{equation}
    \lozenge ( \phi \vee \psi) \; \leftrightarrow \; \lozenge \phi \vee \lozenge \psi
\end{equation}
\begin{equation}
  \Box \phi \vee \Box \psi \; \rightarrow \; \Box ( \phi \vee \psi )
\end{equation}
\begin{equation}
  \lozenge ( \phi \wedge \psi) \; \rightarrow \; \lozenge \phi \wedge
  \lozenge \psi
\end{equation}
\begin{equation}
    \Box ( \phi \rightarrow \psi ) \; \rightarrow \; ( \Box \phi
    \rightarrow \Box \psi )
\end{equation}
\begin{equation}
    ( \lozenge \phi \rightarrow \lozenge \psi ) \; \rightarrow \;
    \lozenge ( \phi \rightarrow \psi)
\end{equation}

\bigskip

\begin{definition}
\end{definition}
\emph{S4 formal system:}

the formal system obtained adding to the \emph{T formal system}
the following axioms:
\begin{equation}
    \Box \Box \phi \; \leftrightarrow \; \Box \phi
\end{equation}
\begin{equation}
    \lozenge \lozenge \phi \; \leftrightarrow \; \lozenge \phi
\end{equation}

\bigskip

\begin{definition} \label{def:S5}
\end{definition}
\emph{S5 formal system:}

the formal system obtained adding to the \emph{S4 formal system}
the following axiom:
\begin{equation}
    \lozenge \Box \phi \; \rightarrow \; \Box \phi
\end{equation}

\newpage
\section{G\"{o}del's ontological proof}

Kurt G\"{o}del formalized (Leibniz's elaboration of Descartes'
elaboration of) the ontological proof of the existence of God
furnished by \emph{Anselmo da Aosta} in his \emph{Proslogion}   as
a theorem of a suitable formal system of Modal Logic though he
didn't published his result since, according to Morgenstern, he
was afraid that his purely logical investigation could be seen as
a religious affair.

 G\"{o}del ontological proof has survived his death as many of his
 unpublished works \cite{Godel-95} and owing to the fact that in 1970 G\"{o}del showed his proof to Dana Scott and allowed him to make a copy of his handwriting pages photocopies
 of which began to circulate in the early eighties, making his first public appearance in
 \cite{Sobel-87}.

The idea of Anselm's ontological proof is that, defined God as the
entity having every perfection (i.e. every ethically and
aesthetically positive property), it must exist since otherwise
one could think a more perfect being having also the perfection of
existing.

Descartes remarked how the essence of such an ontological proof
lies in the fact that God is defined as en entity possessing the
property of \emph{necessary existence}, i.e. the property that if
it is possible than it is necessary.

Leibniz remarked how Descartes' reformulation of Anselm's argument
as a proof of the fact that, by definition, God possesses the
property of necessary existence, had to be augmented with the
proof that the existence of God is possible.

In the framework of the \emph{S5 formal system} let us introduce a
\emph{positivity predicate} $P(\phi) $ and let us assume the
following:
\begin{axiom} \label{ax:first}
\end{axiom}
\begin{equation}
    P( \neg \phi ) \; \leftrightarrow \; \neg P( \phi )
\end{equation}
\begin{axiom} \label{ax:second}
\end{axiom}
\begin{equation}
    P(\phi) \wedge \forall x [ \phi(x) \rightarrow \psi(x) ] \; \rightarrow
    \; P(\psi)
\end{equation}

Then:
\begin{theorem} \label{th:first}
\end{theorem}
\begin{equation}
    P(\phi) \; \rightarrow \; \lozenge \exists x \: \phi(x)
\end{equation}
\begin{proof}
Let as assume the  \emph{ad absurdum hypothesis} $ P(\phi) \wedge
 \neg \lozenge \exists x \phi(x)$.

By the Duns Scoto's principle \emph{ex absurdo quodlibet
sequitur}:
\begin{equation}\label{eq:first}
  \neg \lozenge \exists x \: \phi(x\phi) \; \rightarrow \; \forall \psi \Box
  \forall x [ \phi(x) \rightarrow \psi(x) ]
\end{equation}
Choosing in particular $ \psi := \neg \phi $ the equation
\ref{eq:first} becomes:
\begin{equation}
  \neg \lozenge \exists x \: \phi(x\phi) \; \rightarrow \;
  \forall x [ \phi(x) \rightarrow  \neg \phi(x) ]
\end{equation}
By the axiom \ref{ax:second} it follows that:
\begin{equation}\label{eq:second}
    P( \phi) \wedge \Box \forall x [ \phi(x) \rightarrow ( \neg
    \phi) x \; \rightarrow \; P ( \neg \phi )
\end{equation}
By the axiom \ref{ax:first}:
\begin{equation}
    P ( \neg \phi ) \rightarrow \neg P( \phi )
\end{equation}
so that we obtain that $ \neg P( \phi ) $ contradicting the
hypothesis
\end{proof}

Let us now define the predicate of \emph{being God-like} as the
condition of having all the positive properties:

\begin{definition} \label{def:God-like}
\end{definition}
\begin{equation}
    G(x) \; := \; \forall \phi [ P( \phi ) \rightarrow \phi (x) ]
\end{equation}

Let us now assume  that to be God-like is positive:
\begin{axiom} \label{ax:third}
\end{axiom}
\begin{equation}
    P(G)
\end{equation}

Let us assume furthermore that to be positive cannot be a
contingent property:
\begin{axiom} \label{ax:fourth}
\end{axiom}
\begin{equation}
    P( \phi ) \; \rightarrow \; \Box P( \phi )
\end{equation}

Let us now define the essence of a entity as a property of that
entity implying any other property of its:
\begin{definition} \label{def:essence}
\end{definition}
\begin{equation}
    \phi \: ess \, x \; := \; \phi(x) \wedge \forall \psi \{
    \psi(x) \rightarrow \Box \forall y [ \phi(y) \rightarrow
    \psi(y) ] \}
\end{equation}

Let us now prove that if an entity is God-like to be God-like is
its essence:

\begin{theorem} \label{th:second}
\end{theorem}
\begin{equation}
    G(x) \; \rightarrow \; G \: ess \, x
\end{equation}
\begin{proof}

Applying the following:
\begin{lemma} \label{lem:first}
\end{lemma}
\begin{equation}
 \psi (x) \; \rightarrow \; P( \psi )
\end{equation}
it follows by the axiom \ref{ax:fourth} that:
\begin{equation}
    P( \psi ) \; \rightarrow \; \Box P( \psi )
\end{equation}
By the following:
\begin{lemma} \label{lem:second}
\end{lemma}
\begin{equation}
    \forall x \{ G(x) \leftrightarrow \forall \phi [ P( \phi ) \rightarrow \phi ( x)]
    \} \; \rightarrow \; \Box \{ P( \psi) \rightarrow \forall x
    [G(x) \rightarrow \psi (x) ] \}
\end{equation}
and by the definition \ref{def:God-like} and the definition
\ref{def:essence} it follows that:
\begin{equation}
    \Box \forall x \{ G(x) \leftrightarrow \forall \phi [ P(\phi) \rightarrow \psi(x) ]   \}
\end{equation}
By the lemma \ref{lem:second}:
\begin{equation}
    \Box \{ P( \phi ) \rightarrow \forall x [ G(x) \rightarrow
    \psi(x) ] \}
\end{equation}
and hence, by the definition \ref{def:S5}:
\begin{equation}
    \Box P( \psi ) \; \rightarrow \; \Box \forall x [ G(x)
    \rightarrow \psi (x) ]
\end{equation}
By \emph{modus ponens}:
\begin{equation}
    \Box \forall x [ G(x) \rightarrow \psi (x) ]
\end{equation}
Hence:
\begin{equation}
    G(x) \wedge \forall \psi \{ \psi (x) \rightarrow \Box \forall y
    [ G(y) \rightarrow \psi(y) ] \}
\end{equation}
and hence $ G \: ess \, x $
\end{proof}

Let us now introduce the notion of \emph{necessary existence}:

\begin{definition} \label{def:necessary existence}
\end{definition}
\begin{equation}
    NE( x) \; := \; \forall \phi [ \phi \; ess \, x \; \rightarrow
    \, \Box \exists \phi(x) ]
\end{equation}

Assuming that having the property of necessary existence is
positive:
\begin{axiom} \label{ax:fifth}
\end{axiom}
\begin{equation}
    P(NE)
\end{equation}
we are at last able to prove the following::
\begin{theorem} \label{th:Godel's ontological theorem}
\end{theorem}
\emph{G\"{o}del's Ontological Theorem}
\begin{equation}
    \Box \exists x G(x)
\end{equation}
\begin{proof}
\begin{lemma} \label{lem:third}
\end{lemma}
\begin{equation}
    G(x) \; \rightarrow \; NE(x) \wedge G \: ess \, x
\end{equation}
\begin{lemma} \label{lem:fourth}
\end{lemma}
\begin{equation}
    \exists x G(x) \; \rightarrow \; \Box  \exists x G(x)
\end{equation}
\begin{lemma} \label{lem:fifth}
\end{lemma}
\begin{equation}
    \lozenge \exists x G(x) \; \rightarrow \; \lozenge \Box \exists x G(x)
\end{equation}
\begin{lemma} \label{lem:sixth}
\end{lemma}
\begin{equation}
    \lozenge \exists x G(x) \; \rightarrow \; \Box \exists x G(x)
\end{equation}
\begin{lemma} \label{lem:seventh}
\end{lemma}
\begin{equation}
   \lozenge \exists x G(x)
\end{equation}
\end{proof}

\bigskip

In order to discuss  briefly the literature concerning G\"{o}del's
ontological proof, let us introduce the following:

\begin{definition} \label{def:Godel's ontological formal system}
\end{definition}
\emph{G\"{o}del's ontological formal system:}

\begin{center}
 the formal system $ \mathcal{O} $ consisting in the  \emph{S5  formal system} of the definition \ref{def:S5} augmented
 with the axiom \ref{ax:first} ,the axiom \ref{ax:second}, the axiom
 \ref{ax:third}, the axiom \ref{ax:fourth} and the axiom \ref{ax:fifth}
\end{center}

A disturbing property of G\"{o}del's ontological formal system has
been remarked by Jordan Howard Sobel \cite{Sobel-04}:

\begin{theorem} \label{th:Sobel's theorem}
\end{theorem}
\emph{Sobel's Theorem on modal collapse:}

\begin{hypothesis}
\end{hypothesis}
\begin{equation*}
  \mathcal{O}
\end{equation*}

\begin{thesis}
\end{thesis}
\begin{equation}
    \lozenge \phi \; \Rightarrow \; \Box \phi
\end{equation}

Modifications of G\"{o}del's ontological proof aimed to bypass the
problems related to the theorem \ref{th:Sobel's theorem} have been
proposed by C. Anthony Anderson
\cite{Anthony-Anderson-Gettings-96}.

\newpage
\section{Temporalized ontological proof}
Let us now introduce the temporal operators of Arthur Prior's
Temporal Logic \cite{Prior-03}:
\begin{itemize}
    \item $ \forall_{-} \; := $ "it has been always true that"
    \item  $ \forall_{+} \; := $ "it will be always true that"
    \item $ \exists_{-} \; := \; $ "it has been  true that"
    \item $ \exists_{+} \; := \; $ "it will be   true that"
\end{itemize}

Let us then introduce the following:
\begin{definition}
\end{definition}
\emph{temporal logic:}

the modal logic obtained from the propositional logic introducing
the temporal operators $ \forall_{-} , \forall_{+} , \exists_{-} ,
\exists_{+} $ satisfying the axioms:
\begin{equation}
  \forall_{-} \phi \; \rightarrow \; \exists_{-} \phi
\end{equation}
\begin{equation}
  \forall_{+} \phi \; \rightarrow \; \exists_{+} \phi
\end{equation}

Let us now introduce the following:
\begin{definition} \label{def:temporalization operator}
\end{definition}
\emph{temporalization operator:}
\begin{equation}
    \forall \phi \; \stackrel{\mathcal{T}}{\rightarrow} \;
    \forall_{-} \phi
    \wedge  \forall_{+} \phi
\end{equation}
\begin{equation}
    \exists  \phi\; \stackrel{\mathcal{T}}{\rightarrow} \; \exists_{-} \phi   \wedge
    \exists_{+} \phi
\end{equation}

Let us now introduce the following:

\begin{definition}
\end{definition}
\emph{time reversal operator:}

the operator acting on the temporal labels of the temporal
quantificators in the following way:
\begin{equation}
    \forall_{-} \phi \; \stackrel{T}{\mapsto} \; \forall_{+} \phi
\end{equation}
\begin{equation}
    \exists_{-} \phi \; \stackrel{T}{\mapsto} \; \exists_{+} \phi
\end{equation}
\begin{equation}
    \forall_{+} \phi \; \stackrel{T}{\mapsto} \; \forall_{-} \phi
\end{equation}
\begin{equation}
    \exists_{+} \phi \; \stackrel{T}{\mapsto} \; \exists_{-} \phi
\end{equation}

By construction:
\begin{theorem}
\end{theorem}
\emph{time-reversal invariance of temporalized formal systems:}
\begin{equation}
    T \mathcal{T} S \; = \;  \mathcal{T} S \; \; \forall S
\end{equation}

In particular let us introduce the following:
\begin{definition} \label{def:termporalized Godel's ontological formal system}
\end{definition}
\emph{temporalized G\"{o}del's ontological formal system:}
\begin{equation}
    \mathcal{O}_{\mathcal{T}} \; := \; \mathcal{T} \mathcal{O}
\end{equation}

\newpage
\section{God's death as a breaking of time-reversal symmetry in the temporalized ontological proof}

The logical problem of \emph{theodicy}, namely the problem of
bypassing the logical incompatibility between the definition of
God as the maximally good entity and the existence of evil in the
world was first raised by Epicurus and has been at the center of
Rational Theology from Leibniz's treatise to recent time.

Taking literally Theodor Wiesengrund Adorno's remark that every
metaphysical notion becomes impotent in front and after Auschwitz,
we can think that it affects the same  concept of perfection and
hence the same Anselm's definition of God as an entity having
every perfection.

This lead to think radically about the title of a celebrated book
by Hans Jonas, i.e. the same concept of God after Auschwitz loses
its meaning.

These philosophical remarks much in the spirit of the contemporary
\emph{theologies of God's death} may be mathematically formalized
as a time-reversal symmetry breaking within the formal system
\ref{def:termporalized Godel's ontological formal system}.

With this regard let us first of all introduce the temporal
analogue of Godel's ontological proof:

\begin{theorem} \label{th:Temporalized ontological proof}
\end{theorem}
\emph{Temporalized ontological proof:}

\begin{hypothesis}
\end{hypothesis}
\begin{equation*}
  \mathcal{O}_{\mathcal{T}}
\end{equation*}

\begin{thesis}
\end{thesis}
\begin{equation}
    \Box \exists_{-} x G(x)
\end{equation}
\begin{equation}
    \Box \exists_{+} x G(x)
\end{equation}
\begin{proof}
  It is sufficient to combine the definition
  \ref{def:temporalization
  operator} and the theorem \ref{th:Godel's ontological theorem}
\end{proof}

Let us now formalize the breaking of time-reversal symmetry
through the following:

\begin{definition}
\end{definition}
\emph{time-reversal breaking temporalization's operator:}
\begin{equation}
    \forall \phi \; \stackrel{\mathcal{T}_{\mathcal{B}}}{\rightarrow} \;
    ( \forall_{-} \phi
    \wedge  \neg \forall_{+} \phi ) \vee  ( \forall_{+} \phi
    \wedge  \neg \forall_{-} \phi )
\end{equation}
\begin{equation}
    \exists \phi \; \stackrel{\mathcal{T}_{\mathcal{B}}}{\rightarrow} \;
    ( \exists_{-} \phi
    \wedge  \neg \exists_{+} \phi ) \vee  ( \exists_{+} \phi
    \wedge  \neg \exists_{-} \phi )
\end{equation}

Let us now introduce the following:
\begin{definition} \label{def:time-reversal breaking termporalized Godel's ontological formal system}
\end{definition}
\emph{time-reversal breaking temporalized G\"{o}del's ontological
formal system:}
\begin{equation}
    \mathcal{O}_{\mathcal{T}_{\mathcal{B}}} \; := \; \mathcal{T}_{\mathcal{B}} \mathcal{O}
\end{equation}

Then:
\begin{theorem} \label{th:Time-reversal breaking temporalized ontological proof}
\end{theorem}
\emph{Time-reversal breaking temporalized ontological proof:}

\begin{hypothesis}
\end{hypothesis}
\begin{equation*}
  \mathcal{O}_{\mathcal{T}_{B}}
\end{equation*}

\begin{thesis}
\end{thesis}
\begin{equation}
    ( \Box \exists_{-} x G(x) \wedge \Box \neg \exists_{+} x G(x)
    ) \vee ( \Box \exists_{+} x G(x) \wedge \Box \neg \exists_{-} x G(x)
    )
\end{equation}
\begin{proof}
  It is sufficient to combine the definition \ref{def:time-reversal breaking termporalized Godel's ontological formal system} and the theorem \ref{th:Godel's ontological theorem}
\end{proof}

The impossibility of speaking of any kind of perfection after
Auschwitz may be formalized as the following:
\begin{axiom} \label{ax:Impossibility of positivity in the future}
\end{axiom}
\emph{Impossibility of positivity in the future:}
\begin{equation}
    \Box \forall_{+} \neg P (\phi )
\end{equation}

Then:
\begin{theorem}
\end{theorem}
\emph{Theorem about God's death:}

\begin{hypothesis}
\end{hypothesis}
\begin{center}
  $ \mathcal{O}_{\mathcal{T}_{B}} \wedge $  axiom \ref{ax:Impossibility of positivity in the future}
\end{center}
\begin{thesis}
\end{thesis}
\begin{equation}
    \Box \exists_{-} x G(x) \wedge \Box \neg \exists_{+} x G(x)
\end{equation}
\begin{proof}
The thesis trivially follows combining the theorem
\ref{th:Time-reversal breaking temporalized ontological proof} and
the axiom  \ref{ax:Impossibility of positivity in the future}
\end{proof}

\end{document}